\documentclass[a4paper,12pt]{article}
\usepackage{amsmath,amsfonts,amssymb,amsthm}
\newtheorem{theo}{Th\'eor\`eme}

\newtheorem{lem}{Lemme}

\begin{document}
\title{Une identit\'e remarquable en th\'eorie des partitions}
\author{Alain Lascoux\\{\small Centre National de la Recherche Scientifique}\\
{\small Institut Gaspard Monge, Universit\'e de Marne-la-Vall\'ee}\\
{\small 77454 Marne-la-Vall\'ee Cedex, France}\\
{\small e-mail: Alain.Lascoux @ univ-mlv.fr}\\
Michel Lassalle\\{\small Centre National de la Recherche Scientifique}\\
{\small Ecole Polytechnique}\\
{\small 91128 Palaiseau, France}\\
{\small e-mail: lassalle @ chercheur.com}}
\date{}
\maketitle

\begin{abstract}
We prove an identity about partitions, previously conjectured in the study of shifted Jack polynomials. The proof given is using 
$\lambda$-ring techniques. It would be interesting to obtain a bijective proof.
\end{abstract}

\section{Notations}

Nous d\'emontrons dans cet article une conjecture 
pr\'esent\'ee  dans un pr\'ec\'edent travail ~\cite{La1}. Il s'agit d'une identit\'e  
qui se rencontre dans l'\'etude des polyn\^omes 
``sym\'etriques d\'ecal\'es'' ~\cite{La3,La4}, o\`u elle 
permet le d\'eveloppement explicite de certains ``polyn\^omes 
de Jack d\'ecal\'es'', notamment ceux assoc\'es aux partitions lignes 
et colonnes.

Cette identit\'e se formule de mani\`ere extr\^emement simple dans 
le cadre de la th\'eorie  classique des partitions. 
Cependant il nous a sembl\'e que sa preuve ne s'obtient commod\'ement 
qu'en utilisant la structure (\'el\'ementaire) de $\lambda$-anneau
de l'anneau des polyn\^omes. 

Une partition $\lambda$ est une suite d\'ecroissante finie d'entiers positifs. On dit 
que le nombre
 $n$ d'entiers non nuls est la longueur de $\lambda$. On note
$\lambda  = ( {\lambda }_{1},...,{\lambda }_{n})$ 
et $n = l(\lambda)$. On dit que 
$\left|{\lambda }\right| = \sum\limits_{i = 1}^{n} {\lambda }_{i}$
est le poids de $\lambda$, et pour tout entier $i\geq1$ que 
${m}_{i} (\lambda)  = \textrm{card} \{j: {\lambda }_{j}  = i\}$
est la multiplicit\'e de $i$ dans $\lambda$. On identifie $\lambda$ \`a son diagramme 
de Ferrers 
$\{ (i,j) : 1 \le i \ \le l(\lambda), 1 \ \le j \ \le {\lambda }_{i} \}$. 
On pose 
\[{z}_{\lambda }  = \prod\limits_{i \ge  1}^{} {i}^{{m}_{i}(\lambda)} {m}_{i}(\lambda) 
!  .\]
	 
La g\'en\'eralisation suivante du coefficient 
binomial classique a \'et\'e introduite dans~\cite{La2}. Soient $\lambda$ une partition et $r$ un entier $\geq 1$. On note 
$\genfrac{\langle}{\rangle}{0pt}{}{\lambda}{r}$ 
le nombre de fa\c cons dont on peut choisir r points 
dans le diagramme de $\lambda$ de telle sorte que 
\textit{au moins un point soit choisi sur chaque ligne de $\lambda$}.
	
Les coefficients binomiaux g\'en\'eralis\'es  
$\genfrac{\langle}{\rangle}{0pt}{}{\lambda}{r}$ poss\`edent la fonction 
g\'en\'eratrice suivante
\[\sum_{r \ge  1} \genfrac{\langle}{\rangle}{0pt}{}{\lambda}{r}\,{q}^{r} = 
\prod_{i = 1}^{l(\lambda)} \left({{(1 + q)}^{{\lambda}_{i}} - 1}\right)
 = \prod_{i \ge  1}^{} {\left({{(1 + q)}^{i} - 1}\right)}^{{{m}_{i}(\lambda)}}.
\]

Soient $z$ une ind\'etermin\'ee et $n$ un entier $\geq1$. On note d\'esormais 
\[{(z)}_{n }  =  z (z +1) ... (z + n-1) \quad , \quad {[z]}_{n }  =  z (z -1) ... (z -n +1)\]	
les factorielles ``ascendante'' et ``descendante'' classiques. On pose 
\[\binom{z}{n}  =  {\frac{{[z]}_{n}}{n!}}  .\]

Soit $X=\{X_1,X_2,X_3,\ldots\}$ une famille (infinie) d'ind\'etermin\'ees ind\'ependantes. 
Pour tous entiers $j,k \ge 0$ on pose
\begin{equation}
P_{jk}(X)=\sum_{|\mu|=j} 
\frac{\displaystyle{\genfrac{\langle}{\rangle}{0pt}{}{\mu}{k}}}{z_\mu}\prod_{i\ge 
1}{X_i}^{m_{i}(\mu)}.
\end{equation}

Comme on a $\genfrac{\langle}{\rangle}{0pt}{}{\mu}{k}  = 0$ si  
$ k < l(\mu )$, la sommation est limit\'ee aux partitions $\mu$ telles 
que $l(\mu) \le k$. Il en r\'esulte que $P_{jk}(X)$ est un polyn\^ome de degr\'e $k$. 
Comme on a $\genfrac{\langle}{\rangle}{0pt}{}{\mu}{k}  = 0$ si  
$k > \left|{\mu }\right|$, on a $P_{jk}(X)=0$ pour tout $k > j$. On pose par 
convention $P_{00}(X)=1$.

On a par exemple facilement 
\[P_{j1}(X)= X_{j},\]
\[P_{j2}(X)=\frac{1}{2}(j-1)X_{j}+\frac{1}{2}
\sum_{\begin{subarray}{1}j_{1}+j_{2}=j\\j_{1},j_{2} \ge 1\end{subarray}}X_{j_{1}}X_{j_{2}}.\]

\section{Notre r\'esultat}

Le but de cet article est de d\'emontrer la conjecture suivante, que 
l'un de nous a formul\'ee  dans un pr\'ec\'edent travail (~\cite{La1}, Conjecture 2). 
Cette conjecture explicite un d\'eveloppement en s\'erie formelle. 

\begin{theo}
Soient $z,u$ et $X=\{X_1,X_2,X_3,\ldots\}$
des ind\'etermin\'ees ind\'ependantes. Pour tous entiers $n,r \ge 1$ on a
\begin{multline*}
\sum_{\left|{\mu }\right| = n} (-1)^{r-l(\mu)}
\frac{\displaystyle{\genfrac{\langle}{\rangle}{0pt}{}{\mu}{r}}}{z_{\mu}}
\prod_{i \ge 1}
 {\left(z+\sum_{k \ge 1}u^k\frac{{(i)}_{k}}{k!}\,X_{k} 
 \right)}^{m_{i}(\mu)} =\\
\sum_{j \ge 0} u^j \binom{n+j-1}{n-r}\left( \sum_{k=0}^{min(r,j)} 
\binom{z-j}{r-k} P_{jk}(X) 
\right).
\end{multline*}
\end{theo}

Cette conjecture est triviale pour $r > n$ car on a alors
$\genfrac{\langle}{\rangle}{0pt}{}{\mu}{r}  = 0$. Pour $r=n$ on obtient le r\'esultat 
suivant.
\begin{theo}
Soient $z,u$ et $X=\{X_1,X_2,X_3,\ldots\}$
des ind\'etermin\'ees ind\'ependantes. Pour tout entier $n \ge 1$ on a
\begin{multline*}
\sum_{\left|{\mu }\right| = n} \frac{(-1)^{n-l(\mu)}}{z_{\mu}}\prod_{i \ge 1}
 {\left(z+\sum_{k \ge 1}u^k\frac{{(i)}_{k}}{k!}\,X_{k} 
 \right)}^{m_{i}(\mu)} =\\
\sum_{j \ge 0} u^j\left( \sum_{k=0}^{min(n,j)} \binom{z-j}{n-k} 
P_{jk}(X) 
\right).
\end{multline*}
\end{theo}

Le Th\'eor\`eme 2 avait \'et\'e auparavant conjectur\'e dans ~\cite{La2} 
(Conjecture 4, page 465). Pour $X=0$ le Th\'eor\`eme 1 redonne le Th\'eor\`eme 1' de 
~\cite{La2} (page 462).

\section{Fonctions sym\'etriques}

Nous donnons d'abord ici les notations dont nous aurons besoin \`a 
propos de l'alg\`ebre $\mathbf{Sym}$ des fonctions sym\'etriques, 
consid\'er\'ee d'un point de vue formel.

Soit $A=\{a_1,a_2,a_3,\ldots\}$ un ensemble de variables, qui peut \^etre infini 
(nous dirons que $A$ est un alphabet). On introduit les 
fonctions g\'en\'eratrices
\[ \lambda_t(A) = \prod_{a\in A} (1 +ta) \quad , \quad 
 \sigma_t(A)=  \prod_{a\in A}  \frac{1}{1-ta}
\quad , \quad \Psi_t(a)= \sum_{a\in A} \frac{a}{1-ta} \]
dont le d\'eveloppement d\'efinit les fonctions sym\'etriques \'el\'ementaires
$\Lambda^i(A)$, les fonctions compl\`etes 
$S^i(A)$ et les sommes de puissances $\psi^i(A)$~:

\[ \lambda_t(A) =\sum_{i\geq0} t^i\, \Lambda^i(A)  \quad  , \quad 
 \sigma_t(A) = \sum_{i\geq0} t^i\, S^i(A)
\quad  , \quad \Psi_t(A)=\sum_{i\geq1} t^{i-1}\psi^i(A) .\]

Lorsque l'alphabet $A$ est infini, chacun de ces trois ensembles de fonctions forme une base alg\'ebrique de 
$\mathbf{Sym}[A]$, l'alg\`ebre des fonctions sym\'etriques sur $A$ (c'est-\`a-dire 
que ses \'el\'ements sont alg\'ebriquement ind\'ependants). 

On peut donc d\'efinir l'alg\`ebre $\mathbf{Sym}$ des fonctions sym\'etriques, 
sans r\'ef\'erence \`a l'alphabet $A$, comme l'alg\`ebre sur 
$\mathbf{Q}$ engendr\'ee par les fonctions $\Lambda^i$, $S^i$  
ou $\psi^i$.  

Pour toute partition $\mu=(\mu_{i},1 \le i \le 
l(\mu))=(i^{m_{i}(\mu)},i \ge 1)$, on d\'efinit les fonctions $\Lambda^\mu$, $S^\mu$ 
ou $\psi^\mu$ en posant
\[f^{\mu}= \prod_{i=1}^{l(\mu)}f^{\mu_{i}}=\prod_{k\geq1}(f^k)^{m_{k}(\mu)},\]
o\`u $f^i$ d\'esigne respectivement $\Lambda^i$, $S^i$  ou $\psi^i$. 
Les fonctions $\Lambda^\mu$, $S^\mu$, $\psi^\mu$ forment une base 
lin\'eaire de l'alg\`ebre $\mathbf{Sym}$.

On a la formule de Cauchy
\[\Lambda^{i} = \sum_{\left|{\mu }\right| = i}(-1)^{i-l(\mu)} \frac{\psi^{\mu}}{z_{\mu}}\]
ou encore
\[S^{i} = \sum_{\left|{\mu }\right| = i} \frac{\psi^{\mu}}{z_{\mu}}.\]

Pour toute partition $\mu$, on peut d\'efinir les fonctions sym\'etriques monomiales 
$\psi_\mu$ et les fonctions de Schur $S_\mu$, qui forment \'egalement une base 
lin\'eaire de l'alg\`ebre $\mathbf{Sym}$.

Les bases $\Lambda^\mu$, $S^\mu$, $\psi^\mu$, $\psi_\mu$ ou $S_\mu$ 
sont not\'ees respectivement $e_{\mu}$, $h_{\mu}$, 
$p_{\mu}$, $m_{\mu}$ ou $s_{\mu}$ dans la litt\'erature, notamment dans ~\cite{Ma}. 
Les notations utilis\'ees ici sont celles de ~\cite{LS}, qui sont plus 
adapt\'ees aux $\lambda$-anneaux. 

\section{L'anneau des polyn\^omes comme $\lambda$-anneau}

Nous allons d\'emontrer le Th\'eor\`eme 1 en utilisant le fait que 
l'anneau des polyn\^omes poss\`ede une structure de $\lambda$-anneau. 

Un $\lambda$-anneau est un anneau commutatif avec unit\'e muni 
d'op\'erateurs qui v\'erifient certains axiomes. Nous renvoyons le lecteur \`a  
~\cite{K} pour la th\'eorie g\'en\'erale, et au chapitre 2  de ~\cite{Pr} 
pour leur application \`a l'analyse multivari\'ee.

Nous n'utiliserons cette th\'eorie que dans le cadre \'el\'ementaire 
suivant. Soit $A=\{a_1,a_2,a_3,\ldots\}$ un alphabet quelconque. On consid\`ere l'anneau  $\mathbf{R}[A]$ des 
polyn\^omes en $A$ \`a coefficients r\'eels. La structure de  
$\lambda$-anneau de $\mathbf{R}[A]$ consiste \`a 
d\'efinir une action de $\mathbf{Sym}$ sur $\mathbf{R}[A]$.

\subsection{Action de $\mathbf{Sym}$}

Les fonctions $\psi^i$ formant un syst\`eme de g\'en\'erateurs alg\'ebriques 
de $\mathbf{Sym}$, \'ecrivant tout polyn\^ome sous la forme $\sum_{c,u} c u$,
avec $c$ constante r\'eelle et $u$ un mon\^ome en 
$\{a_{1},a_{2},a_{3}\ldots\}$, on d\'efinit une action de $\mathbf{Sym}$ 
sur $\mathbf{R}[A]$, not\'ee $[.]$, 
en posant 
\[ \psi^{i} [\sum_{c,u} c u]=\sum_{c,u} 
cu^{i}.\]

Pour tous polyn\^omes $P,Q \in \mathbf{R}[A]$ on en d\'eduit 
imm\'ediatement $\psi^{i} [PQ]=\psi^{i} [P] \psi^{i} [Q]$ et
$\psi^{\mu} [PQ]=\psi^{\mu} [P] \psi^{\mu} [Q]$.

L'action ainsi d\'efinie s'\'etend \`a tout \'el\'ement de $\mathbf{Sym}$. Ainsi on a
\[ \lambda_t [\sum_{c,u} c u]= \prod_{c,u} (1+tu)^c \quad , \quad
\sigma_t [\sum_{c,u} c u]= \prod_{c,u} (1-tu)^{-c}.\]
On en d\'eduit $S^{i} [P] =(-1)^i \Lambda^{i} [-P]$.

On notera le comportement diff\'erent des 
constantes $c\in \mathbf{R}$ et des mon\^omes $u$ :
\begin{equation}
\begin{split}   
\psi_i [c] = c \quad , \quad  S^i [c] =\frac{(c)_{i}}{i!} \quad , \quad  
\Lambda^i [c] = \frac{[c]_{i}}{i!} \quad \\
\psi_i [u] = u^i = S^i [u]  \quad , \quad  
\Lambda^{i} [u] = 0,\, i>1 \quad , \quad \Lambda^{1} [u] = u.
\end{split}
\end{equation}

Il est plus correct de caract\'eriser les "mon\^omes" $u$ comme \'el\'ements {\it  
de rang 1} (i.e. les $u\neq 0,1$ tels que $\Lambda^i[u] = 0 \ \forall i>1$),
et les "constantes" $c\in \mathbf{R}$ comme les \'el\'ements invariants par les $\psi_i$
(on dira aussi \'el\'ement {\it de type binomial}).

Lorsqu'on utilise la th\'eorie des $\lambda$-anneaux pour d\'emontrer 
une identit\'e alg\'ebrique, il est donc \textit{toujours n\'ecessaire de 
pr\'eciser le statut de chaque \'el\'ement}. En particulier nous 
aurons \`a employer des ind\'etermin\'ees de rang 1, et d'autres de type binomial.

\subsection{Extension aux s\'eries formelles}

On remarquera que si $a_{1},a_{2},\ldots,a_{n}$ sont des \'el\'ements de 
rang un, alors  
\[\psi^{i} [a_{1}+a_{2}+\ldots+a_{n}]=a_{1}^i+a_{2}^i+\ldots+a_{n}^i\]
est la valeur de la $i$-\`eme somme de puissance $\psi^{i} 
(a_1,a_2,\ldots,a_{n})$. 

Dans la suite pour tout alphabet 
$A=\{a_1,a_2,a_3,\ldots\}$, on notera $A^\diamondsuit =\sum_{i}a_{i}$ la 
somme de ses \'el\'ements. Lorsque A est form\'e d'\'el\'ements de rang 
1, on a ainsi pour toute fonction sym\'etrique $f$,
\begin{equation}
f[A^\diamondsuit]=f(A).
\end{equation}

En particulier si $q$ est de rang 1, on a 
\[\psi^{i} (1,q,q^2,q^3,\ldots,q^{n-1})= \psi^{i} [\sum_{k=0}^{n-1} q^k].\]
Il est naturel de vouloir \'ecrire
\[\sum_{k=0}^{n-1} q^k=\frac{1-q^n}{1-q},\] 
et d'\'etendre ainsi l'action de $\mathbf{Sym}$ aux fonctions rationnelles.  
Il est \'egalement naturel de consid\'erer un alphabet infini 
$(1,q,q^2,q^3,\ldots)$, de vouloir sommer la s\'erie 
\[\sum_{k\geq0}q^k=\frac{1}{1-q},\] 
et d'\'etendre ainsi l'action de $\mathbf{Sym}$ aux s\'eries formelles \`a 
coefficients r\'eels.

Pour cela on pose
\[ \psi_i \left( \frac{\sum c u}{\sum d v}\right) = 
  \frac{\sum c u^i}{\sum d v^i}  \ , \]
avec $c,d$ constantes r\'eelles et $u,v$ des mon\^omes en 
$(a_{1},a_{2},a_{3}\ldots)$.

L'action ainsi d\'efinie s'\'etend \`a tout \'el\'ement de $\mathbf{Sym}$. 
On munit ainsi l'anneau des s\'eries formelles \`a 
coefficients r\'eels d'une structure de $\lambda$-anneau.
On a par exemple
\[\lambda_t [\frac{1}{1-q}]= \prod_{i\geq1} (1+tq^i)={(-t;q)}_{\infty}\]
\[\sigma_t [\frac{1}{1-q}]= \prod_{i\geq1} 
\frac{1}{1-tq^i}=\frac{1}{{(t;q)}_{\infty}},\]
ce qui fait appara\^itre des quantit\'es bien connues en $q$-calcul. 

\subsection{Formulaire}

Les relations fondamentales suivantes sont des 
cons\'equences directes des relations (2). Certaines nous seront 
n\'ecessaires. Pour tous $P,Q$ on a d'abord
\begin{equation}
\begin{split}
S^i [P+Q]&= \sum_{j=0}^{i} S^{i-j} [P] S^j [Q]\\
\Lambda^{i} [P+Q]&= \sum_{j=0}^{i} \Lambda^{i-j} [P] \Lambda^j [Q] ,
\end{split}
\end{equation}
ou de mani\`ere \'equivalente :
\begin{equation}
\begin{split}
\sigma_t[P+Q]=\sigma_t[P]\,\sigma_t[Q] \\
\lambda_t[P+Q]=\lambda_t[P]\,\lambda_t[Q].
\end{split}
\end{equation}

Pour tous $P,Q$ on a d'autre part
\begin{equation}
\begin{split}
S^i [PQ]&= \sum_{\left|{\mu }\right| = i} 
\frac{1}{z_{\mu}}\psi^{\mu} [P] \psi^{\mu} [Q]\\
&=\sum_{\left|{\mu }\right| = i} \psi_{\mu} [P] S^{\mu} [Q]\\
&=\sum_{\left|{\mu }\right| = i} S_{\mu} [P] S_{\mu} [Q],
\end{split}
\end{equation}
ou de mani\`ere \'equivalente :
\begin{equation}
\begin{split}
\Lambda^{i} [PQ]&= \sum_{\left|{\mu }\right| = i} 
\frac{(-1)^{i-l(\mu)}}{z_{\mu}} \psi^{\mu} [P] \psi^{\mu} [Q]\\
&=\sum_{\left|{\mu }\right| = i} \psi_{\mu} [P] \Lambda^{\mu} [Q]\\
&=\sum_{\left|{\mu }\right| = i} S_{\mu} [P] S_{\mu'} [Q],
\end{split}
\end{equation}
o\`u $\mu'$ d\'esigne la partition transpos\'ee de $\mu$.

Si $P$ est de rang 1 et $Q$ arbitraire, on a
\[\Lambda^{i} [PQ] =P^i \Lambda^{i}[Q].\]
Ainsi lorsque $P$ et $Q$ sont de rang 1, on a
\begin{equation}
\lambda_t[PQ]=1+tPQ.
\end{equation}

\section{D\'emonstration du Th\'eor\`eme 1}    

\subsection{Pr\'eliminaires}
\begin{lem}
Soit $q'$ un \'el\'ement de rang 1. Si on pose $q=q'-1$, on a  \[ \psi^{\mu} [q]=
\sum_{k \ge  1} \genfrac{\langle}{\rangle}{0pt}{}{\mu}{k}\,{q}^{k}.\]
\end{lem}
\begin{proof}[Preuve]
On a
\[ \psi^{i} [q]=\psi^{i} [q'-1]=(q')^i-1={(1+q)}^i-1.\]
On applique la fonction g\'en\'eratrice des entiers 
$\genfrac{\langle}{\rangle}{0pt}{}{\mu}{k}$.
\end{proof}

A l'aide des relations (2) et (4) on obtient facilement
\begin{equation}
\Lambda^i [q] = (-1)^{i-1} q \quad , \quad S^i [q] = (1+q)^{i-1} q \quad , 
\quad i\ge1.
\end{equation}
On en d\'eduit
\[\Lambda^{\mu} [q] = (-1)^{\left|{\mu }\right|-l(\mu )} q^{l(\mu )} \quad , \quad
S^{\mu} [q] = (1+q)^{\left|{\mu }\right|-l(\mu )} q^{l(\mu )} .\]

On rappelle la d\'efinition du polyn\^ome $P_{jk}$ introduit en (1) et 
de la fonction sym\'etrique monomiale $\psi_{\mu}$ (somme de tous 
les mon\^omes diff\'erents ayant pour exposant une permutation de $\mu$).

\begin{lem}
Soit $A=\{a_1,a_2,a_3,\ldots\}$ un alphabet (fini ou infini) quelconque. 
Pour tout $i\ge 1$ on pose $X_{i} = \sum_{a \in A}a^i$. 
Alors pour tous entiers $j,k \ge 0$ on a
\[P_{jk}(-X) = (-1)^{k}\sum_{|{\mu}| = j,l(\mu ) = k} \psi_{\mu} (A).\]
\end{lem} 

\begin{proof}[Preuve]
L'\'egalit\'e \`a \'etablir est une identit\'e alg\'ebrique entre 
polyn\^omes en les $a_{i}$. Elle est enti\`eremnt ind\'ependante de la structure de $\lambda$-anneau 
de l'anneau des polyn\^omes. Pour la d\'emontrer dans le cadre de la 
th\'eorie des $\lambda$-anneaux, nous pouvons donc choisir le statut  
de chacune des ind\'etermin\'ees $a_{i}$. 

Nous pouvons par exemple supposer que tous les \'el\'ements de l'alphabet $A$ sont 
de rang 1. Compte-tenu de (3), la relation \`a d\'emontrer devient 
dans ce cas
\begin{equation}
P_{jk}(-X) = (-1)^{k}\sum_{|{\mu}| = j,l(\mu ) = k} \psi_{\mu} 
[A^\diamondsuit].
\end{equation}

Compte-tenu de (3), on a aussi dans ce cas
\[\psi^{i}[A^\diamondsuit] = \psi^{i} (A) = X_{i} \quad , \quad i\ge 1,\]
d'o\`u pour toute partition $\mu$,
\[\psi^{\mu}[A^\diamondsuit] = \psi^{\mu} (A) = \prod_{i\ge1}{X_{i}}^{m_{i}(\mu)}.\]

La formule de Cauchy (7) implique alors 
\begin{equation*}
\begin{split}
\Lambda^{j} [qA^\diamondsuit]&=\sum_{\left|{\mu }\right| = j} 
\frac{(-1)^{j-l(\mu)}}{z_{\mu}} \psi^{\mu} [q] \psi^{\mu} [A^\diamondsuit]\\
&=\sum_{\left|{\mu }\right| = j} 
\frac{(-1)^{j-l(\mu)}}{z_{\mu}}
\left(\sum_{k \ge  1} 
\genfrac{\langle}{\rangle}{0pt}{}{\mu}{k}\,{q}^{k}\right)
 \prod_{i\ge1}{X_{i}}^{m_{i}(\mu)}\\
&= (-1)^{j} \sum_{k\ge1} P_{jk}(-X) q^k.
\end{split}
\end{equation*}
Et d'autre part on a aussi
\begin{equation*}
\begin{split}
\Lambda^{j} [qA^\diamondsuit]&= \sum_{\left|{\mu }\right| = j}
\psi_{\mu} [A^\diamondsuit] \Lambda^{\mu} [q]\\
&=\sum_{\left|{\mu }\right| = j} \psi_{\mu} [A^\diamondsuit] (-1)^{j-l(\mu)} q^{l(\mu)}.
\end{split}
\end{equation*}
On en d\'eduit (10) par comparaison.
\end{proof}

\subsection{M\'ethode}

\textit{Dans toute la suite de cet article}, on consid\`ere un alphabet (fini ou infini) 
$A=\{a_1,a_2,a_3,\ldots\}$. Pour le moment, nous ne faisons aucune hypoth\`ese sur le statut des 
\'el\'ements de $A$. En particulier nous ne supposons pas que les 
$a_{k}$ sont de rang 1. Pour tout $i\ge 1$ on pose 
\[X_{i} = \psi^{i}(A) = \sum_{a \in A}a^i.\]

On consid\`ere quatre \'el\'ements $q',z,t,u$. On suppose que $z$ est de type binomial et que 
$q'=1+q$ est de rang 1.

Pour d\'emontrer l'identit\'e du Th\'eor\`eme 1, on va montrer l'\'egalit\'e 
des fonctions g\'en\'eratrices de ses deux membres. Plus pr\'ecis\'ement on  \'ecrit 
chaque membre de l'identit\'e du Th\'eor\`eme 1 en changeant les $X_i$ en $-X_i$, et 
on somme sur $n$ et $r$ apr\`es avoir multipli\'e 
par ${(-t)}^n $ ${(-q)}^r$.

L'\'egalit\'e \`a d\'emontrer devient
\begin{multline}
\sum_{n \ge r \ge1} \sum_{|\mu| = n} (-1)^{r-l(\mu)} {(-t)}^n {(-q)}^r
\frac{\displaystyle{\genfrac{\langle}{\rangle}{0pt}{}{\mu}{r}}}{z_{\mu}}
\prod_{i \ge 1}
 {\left(z-\sum_{k \ge 1}u^k\frac{{(i)}_{k}}{k!}\,X_{k} 
 \right)}^{m_{i}(\mu)} =\\
\sum_{n \ge r \ge1} \sum_{j \ge 0} {(-t)}^n {(-q)}^r u^j \binom{n+j-1}{n-r}\left( \sum_{k=0}^{min(r,j)} 
\binom{z-j}{r-k} P_{jk}(-X) 
\right).
\end{multline}

\subsection{Membre de droite}
Compte tenu du Lemme 2, le membre de droite de 
(11) s'\'ecrit, en notant $uA=\{ua_1,ua_2,ua_3,\ldots\}$,
\[\sum_{n \ge r \ge1} \sum_{\nu}  {(-t)}^{n} (-q)^r u^{|\nu|} 
\binom{n+|\nu|-1}{n-r} \binom {z-|\nu|}{r-l(\nu)}  (-1)^{l(\nu)} \psi_\nu (A)=\]
\[\sum_{\nu} (-1)^{l(\nu)}  \psi_\nu (uA) \sum_{n \ge r \ge l(\nu)}   
{(-t)}^{n} (-q)^r  \binom{n+|\nu|-1}{n-r} \binom {z-|\nu|}{r-l(\nu)}.  \]

Mais on a la relation suivante, qui est une autre fa\c con 
d'\'ecrire la formule classique du bin\^ome :
\[\sum_{i\ge j}{\binom{i-1}{j-1}t^{i-j}}=\frac{1}{(1-t)^{j}}.\]
On en d\'eduit imm\'ediatement
\[\sum_{n \ge r}  {(-t)}^{n} \binom{n+|\nu|-1}{n-r} = 
\frac{(-t)^r}{(1+t)^{|\nu|+r}} . \]
Le membre de droite de (11) s'\'ecrit donc
\[\sum_{\nu} (-1)^{l(\nu)}  \psi_\nu (uA) \left(\sum_{r \ge l(\nu)}   
 \frac{(qt)^r}{(1+t)^{|\nu|+r}} \binom {z-|\nu|}{r-l(\nu)}\right).  \]
Ce qui peut se reformuler
\[\sum_{\nu} \psi_\nu (uA) \,   
 \frac{(-qt)^{l(\nu)}}{(1+t)^{|\nu|+l(\nu)}} \, 
 \left(\sum_{k\ge0} {\left(\frac{qt}{1+t}\right)}^k \binom {z-|\nu|}{k}\right).  \]

Finalement le membre de droite de (11) s'\'ecrit
\[\sum_{\nu} \psi_\nu (uA) \,  
\frac{(-qt)^{l(\nu)}}{(1+t)^{|\nu|+l(\nu)}} \, 
{\left(1+\frac{qt}{1+t}\right)}^{z-|\nu|}.  \]
Soit encore en posant $y= -qt/(1+t)$,
\begin{equation}
\sum_{\nu} \psi_\nu (uA) \,  
\frac{y^{l(\nu)}}{(1+t)^{|\nu|}} \,  {(1-y)}^{z-|\nu|}.  
\end{equation}

\subsection{Membre de gauche}

Comme on a $X_{k} = \sum_{a \in A}a^k$, la quantit\'e suivante, \'ecrite au membre de 
gauche de (11), devient
\begin{equation*}
\begin{split}
z-\sum_{k \ge 1}u^k\frac{{(i)}_{k}}{k!}\,X_{k}
&=z-\sum_{k \ge 1} u^k \frac{{(i)}_{k}}{k!}\left(\sum_{a \in A}a^k\right) .
\end{split}
\end{equation*}

On introduit alors l'alphabet 
\[A'=\left\{\frac{1}{1-ua_1},\frac{1}{1-ua_2},\frac{1}{1-ua_3},\ldots\right\}=
\left\{\frac{1}{1-ua}, a \in A\right\}.\] 
\textit{Nous faisons d\'esormais l'hypoth\`ese suivante} : chaque \'el\'ement 
$\frac{1}{1-ua}$ est de rang 1. Sous cette hypoth\`ese on a
\begin{equation*}
\begin{split}
\psi^{i} [\sum_{a\in A} \frac{1}{1-ua}]&= \sum_{a\in A} (1-ua)^{-i}\\
&=\sum_{a\in A} \left(\sum_{k \ge 0}\frac{{(i)}_{k}}{k!}u^k a^k\right) 
\\
&=\sum_{a\in A} \left(1+\sum_{k \ge 1}\frac{{(i)}_{k}}{k!}u^k a^k\right).
\end{split}
\end{equation*}

On introduit l'\'el\'ement \[B=z-\sum_{a\in A} \frac{ua}{1-ua}=z+\sum_{a\in 
A} \left(1-\frac{1}{1-ua}\right).\]
On a ainsi 
\[z-\sum_{k \ge 1}u^k\frac{{(i)}_{k}}{k!}\,X_k= \psi^i [B] .\]
Pour toute partition $\mu$, on en d\'eduit
\[\prod_{i \ge 1}
 {\left(z-\sum_{k \ge 1}u^k\frac{{(i)}_{k}}{k!}\,X_{k} 
 \right)}^{m_{i}(\mu)}= \psi^\mu [B] .\] 

Compte-tenu de cette relation, le membre de gauche de (11) s'\'ecrit
\begin{equation*}
\begin{split}
\sum_{n \ge r \ge1} \sum_{|{\mu }| = n} (-1)^{r-l(\mu)} {(-t)}^n {(-q)}^r
\frac{\displaystyle{\genfrac{\langle}{\rangle}{0pt}{}{\mu}{r}}}{z_{\mu}}
\psi^\mu [B]
&= \sum_{n \ge1} t^n \sum_{\left|{\mu }\right| = n} \frac{(-1)^{n-l(\mu)}}{z_{\mu}} 
\left( \sum_{r \ge  1} 
\genfrac{\langle}{\rangle}{0pt}{}{\mu}{r}\,q^r\right )
  \psi^\mu [B] \\
&= \sum_{n \ge1} t^n \sum_{\left|{\mu }\right| = n} \frac{(-1)^{n-l(\mu)}}{z_{\mu}} 
\psi^\mu [q] \ \psi^\mu [B] \\
&= \sum_{n \ge1} t^n  \Lambda^{n} [qB]\\
&=\lambda_t [qB].
\end{split}
\end{equation*}

La d\'emonstration sera termin\'ee en prouvant que le d\'eveloppement 
(12) est exactement la d\'ecomposition 
de $\lambda_t [qB]$ sur la base des fonctions monomiales $\psi_\nu 
(uA)$.

\subsection{D\'eveloppement de $\lambda_t [qB]$.}

On maintient les notations pr\'ec\'edentes en faisant le changement 
de variables $ua \rightarrow a$. On consid\`ere 
un alphabet $A=\{a_1,a_2,a_3,\ldots\}$ et trois 
\'el\'ements $q',z,t$ avec les hypoth\`eses suivantes :

- on suppose que $z$ est de type binomial et que 
$q'=1+q$ est de rang 1, 

- on suppose que pour tout $a\in A$, l'\'el\'ement $a'=\frac{1}{1-a}$ est de rang 1.

On a maintenant 
\[B=z-\sum_{a\in A} \frac{a}{1-a}=z+\sum_{a\in A} \left(1-\frac{1}{1-a}\right).\]
Nous allons d\'emontrer le Th\'eor\`eme 1 sous la forme suivante.

\begin{theo} 
En posant $y= -qt/(1+t)$, on a
\[\lambda_t [qB] =\sum_{\nu} \psi_\nu (A) \,  
\frac{y^{l(\nu)}}{(1+t)^{|\nu|}} \,  {(1-y)}^{z-|\nu|}.  \]
\end{theo}
\begin{proof}[Preuve]
On a d'abord
\[\lambda_t [qB] =\lambda_t [qz+q\sum_{a\in A} (1-a')]\]
Comme $z$ est de type binomial, on a
\[\lambda_t [qz] = (\lambda_t [q])^{z}. \]
Et d'autre part la relation (9) implique
\[\lambda_t [q]=1+q \sum_{i\ge1} (-1)^{i-1}t^i  = 1+ \frac{qt}{1+t}.\]

Compte-tenu de (5) on en d\'eduit
\begin{equation*}
\begin{split}
\lambda_t [qB] &=\lambda_t [qz] \, \lambda_t [q \sum_{a\in A} 
(1-a')]\\
&= {(1-y)}^{z} \, \lambda_t [q\sum_{a\in A} (1-a')]\\
&= {(1-y)}^{z} \, \prod_{a\in A} \lambda_t [q (1-a')].
\end{split}
\end{equation*}
    
Maintenant on a $q(1-a')=(q'-1)(1-a')=q'-1-q'a'+a'$. Les \'el\'ements 
$q'$ et $a'$ \'etant de rang 1, les relations (5) et (8) impliquent
\[\lambda_t [q (1-a')]=\frac{\lambda_t [q']}{\lambda_t [1]} \,
\frac{\lambda_t [a']}{\lambda_t [q'a']}
=\frac{1+tq'}{1+t} \, \frac{1+ta'}{1+tq'a'}
=(1-y) \frac{1+ta'}{1+t(1+q)a'} \,.\]

Finalement on obtient
\begin{equation*}
\begin{split}
\lambda_t [q B] &= {(1-y)}^{z} \prod_{a\in A} (1-y) 
\frac{1+t-a}{1+t(1+q)-a}\\
&={(1-y)}^{z} \prod_{a\in A}  (1-y)
\left(1-\frac{qt}{1+t+qt-a}\right).
\end{split}
\end{equation*}
Posons alors 
\[v=\frac{1}{(1+t)(1-y)}=\frac{1}{1+t+qt}.\]
La relation pr\'ec\'edente devient
\[\lambda_t [q B]= {(1-y)}^{z} \prod_{a\in A}
\left(1+y\frac{va}{1-va}\right).\]

Maintenant on a
\[\prod_{a\in A}
\left(1+y\frac{va}{1-va}\right) = \sum_{N \subset A} \, \prod_{a\in N} y 
\frac{va}{1-va}.\]
Nous allons utiliser la propri\'et\'e suivante, qui se v\'erifie 
facilement :
\[\sum_{\begin{subarray}{1}N \subset A \\ \mathrm{card} N = n\end{subarray}} \, \prod_{a\in N} y \frac{va}{1-va} = {y}^{n} 
\sum_{l(\nu)=n} v^{|\nu|} \psi_\nu(A).\]
Soit encore
\[\lambda_t [q B]= {(1-y)}^{z}  \sum_{\nu}  
\frac{y^{l(\nu)}}{(1+t)^{|\nu|}(1-y)^{|\nu|}} \, \psi_\nu(A).\]
On conclut imm\'ediatement.
\end{proof}

\section{Application}

Les Th\'eor\`emes 1 et 2 peuvent permettre d'obtenir des 
identit\'es remarquables en sp\'ecialisant les ind\'etermin\'ees $X_{i}$ et $z$. 

Nous revenons seulement ici sur les conjectures de 
~\cite{La2}, rencontr\'ees en \'etudiant les 
polyn\^omes sym\'etriques d\'ecal\'es ~\cite {La3,La4}. 
Soit $\alpha$ un nombre r\'eel positif. Pour toute partition $\lambda$ et tout 
entier $k\ge 0$, on note
\[d_{k}(\lambda) = \sum_{(i,j) \in \lambda} {\left(j-1-\frac{i-1}{\alpha}\right)}^k.\]
On introduit la g\'en\'eralisation suivante de la ``factorielle ascendante'':
\[(z)_{\lambda}=\prod_{(i,j) \in \lambda} \left(z+j-1-\frac{i-1}{\alpha}\right).\]

Pour tous entiers $j,k \ge 0$ on pose
\[F_{jk}(\lambda)=P_{jk}(d_{1}(\lambda),d_{2}(\lambda),d_{3}(\lambda),\ldots).\]
C'est-\`a-dire qu'on choisit la sp\'ecialisation suivante  
\[ X_{k}=d_{k}(\lambda) \quad , \quad k \ge 1.\]
En d'autres termes, l'alphabet $A$ tel que $X_{k} = \sum_{a \in A}a^k$ est alors
\[A_{\lambda}=\left\{j-1-\frac{i-1}{\alpha}, (i,j) \in \lambda\right\} .\]

\begin{theo}
Soient $x,y$ deux ind\'etermin\'ees ind\'ependantes. Pour toute 
partition $\lambda$ on a
\[\frac{(y-x)_{\lambda}}{(y)_{\lambda}} =
\sum_{i\ge 0} \sum_{j\ge0} (-1)^{i+j} \frac{x^i}{y^{i+j}}
\left(\sum_{k=0}^{min(i,j)}\binom{|\lambda|-j}{i-k}\,F_{jk}(\lambda) 
\right).\]
\end{theo}

\begin{proof}[Preuve]
On montre comme dans ~\cite{La2} (p. 464) que
\[\frac{(y-x)_{\lambda}}{(y)_{\lambda}} = \sum_{\mu} v^{|\mu|}
 \frac{(-1)^{|\mu|-l(\mu)}}{z_{\mu}}  \prod_{i \ge 1}
 {\left(\sum_{p \ge 0}u^p\frac{{(i)}_{p}}{p!}\,d_{p}(\lambda) 
 \right)}^{m_{i}(\mu)},\]
avec $v=-x/y$ et $u=-1/y$. On \'ecrit le Th\'eor\`eme 2 sp\'ecialis\'e avec $X_{k}=d_{k}(\lambda)$ et $z=d_{0}(\lambda)=|\lambda|$. 
\end{proof}

Il est important de noter que la sommation a lieu sur tout $j\ge 0$ et pas seulement 
sur $|\lambda|-j\ge 0$. Le degr\'e en $x$ du membre de gauche \'etant clairement 
$\le |\lambda|$, on obtient pour tout $i > |\lambda|$,$j\ge 0$,
\[\sum_{k=0}^{min(i,j)}\binom{|\lambda|-j}{i-k}\,F_{jk}(\lambda) = 
0.\]
En effet c'est seulement lorsque l'alphabet $A_{\lambda}$ est infini que les 
ind\'etermin\'ees $d_{k}(\lambda)$ sont ind\'ependantes.

Le cas o\`u $\lambda$ est une partition-ligne $(n)$ correspond au 
d\'eveloppement en s\'erie de la formule classique de Chu-Vandermonde 
~\cite{La2}.

\end{document}